\RequirePackage{ifpdf}
\ifpdf 
\documentclass[pdftex]{sigma}
\else
\documentclass{sigma}
\fi

\usepackage[all]{xy}

\begin{document}


\renewcommand{\thefootnote}{$\star$}

\renewcommand{\PaperNumber}{080}

\FirstPageHeading

\newcommand{\GR}{general relativity}
\newcommand{\MM}{MacDowell--Mansouri}
\newcommand{\tm}{topologically massive}
\newcommand{\tme}{topologically massive electrodynamics}
\newcommand{\tmg}{topologically massive gravity}
\newcommand{\gr}{general relativity}
\newcommand{\cs}{Chern--Simons}
\newcommand{\sgn}{{\rm sgn}}

\newcommand{\bfm}{\boldmath \bf}

\newcommand{\half}{\frac{1}{2}}
\newcommand{\fourth}{\frac{1}{4}}

\newcommand{\xa}{\alpha}
\newcommand{\xb}{\beta}
\newcommand{\xg}{\gamma}
\newcommand{\xd}{\delta}
\newcommand{\xe}{\epsilon}
\newcommand{\om}{\omega}


\newcommand{\R}{{\mathbb R}}
\newcommand{\C}{{\mathbb C}}
\newcommand{\N}{{\mathbb N}}
\newcommand{\Z}{{\mathbb Z}}

\newcommand{\Q}{{\cal Q}}


\def\tensor {{\otimes}}
\def\disunion {{\sqcup}}
\def\bigdisunion {\bigsqcup}
\def\tr {{\rm tr}}


\newcommand{\maps}{\colon}
\def\ker {{\rm ker}}
\def\image {{\rm im}}
\def\ran {{\rm ran}}

\def\stackto #1 { \, {\stackrel{#1}{\longrightarrow}}\, }
\def\stackTo #1 { {\stackrel{#1}{\Longrightarrow}} }
\newcommand{\To}{\Rightarrow}

\def\iso{\cong}

\newcommand{\Vect}{{\rm Vect}}

\newcommand{\frakg}{\mathfrak{g}} %
\newcommand{\frakh}{\mathfrak{h}} %
\newcommand{\frakp}{\mathfrak{p}} %

\newcommand{\G}{\mathcal{G}}

\newcommand{\Ad}{{\rm Ad}}
\newcommand{\ad}{{\rm ad}} 

\renewcommand{\O}{{\rm O}}

\newcommand{\SO}{{\rm SO}}
\newcommand{\so}{\mathfrak{so}}

\newcommand{\ISO}{{\rm ISO}}
\newcommand{\Iso}{\mathfrak{iso}}

\newcommand{\SL}{{\rm SL}}
\newcommand{\SLtilde}{\widetilde{\rm SL}}
\newcommand{\Sl}{\mathfrak{sl}}

\newcommand{\GL}{{\rm GL}}
\newcommand{\gl}{\mathfrak{gl}}

\newcommand{\PSL}{{\rm PSL}}

\newcommand{\U}{{\rm U}}
\renewcommand{\u}{\mathfrak{u}}

\newcommand{\SU}{{\rm SU}}
\newcommand{\su}{\mathfrak{su}}

\newcommand{\Sp}{{\rm Sp}}

\newcommand{\Aut}{{\rm Aut}}
\newcommand{\Der}{{\rm Der}}

\newcommand{\g}{\mathfrak{g}}
\newcommand{\h}{\mathfrak{h}}
\newcommand{\p}{\mathfrak{p}}  

\def\semidir {\ltimes}
\def\subgp {\subset}

\newcommand{\eref}[1]{(\ref{#1})}

\newcommand{\bilin}{\beta} 

\newcommand{\ssp}{{\scriptscriptstyle +}}
\newcommand{\ssm}{{\scriptscriptstyle -}}
\newcommand{\sspm}{{\scriptscriptstyle \pm}}

\renewcommand{\tensor}{\otimes}

\newcommand{\smallhalf}{{\textstyle\frac{1}{2}}}

\newcommand{\Ric}{{\rm Ric}} 

\newcommand{\Spin}{{\rm Spin}}

\newcommand{\beq}{\begin{equation}}
\newcommand{\eeq}{\end{equation}}

\ShortArticleName{Symmetric Space Cartan Connections and Gravity in Three and Four Dimensions}

\ArticleName{Symmetric Space Cartan Connections and Gravity\\ in Three and Four Dimensions\footnote{This paper is a
contribution to the Special Issue ``\'Elie Cartan and Dif\/ferential Geometry''. The
full collection is available at
\href{http://www.emis.de/journals/SIGMA/Cartan.html}{http://www.emis.de/journals/SIGMA/Cartan.html}}}

\Author{Derek K.~WISE}

\AuthorNameForHeading{D.K.~Wise}

\Address{Department of Mathematics, University of California, Davis, CA 95616, USA}
\Email{\href{mailto:derek@math.ucdavis.edu}{derek@math.ucdavis.edu}}
\URLaddress{\url{http://math.ucdavis.edu/~derek/}}

\ArticleDates{Received April 10, 2009, in f\/inal form July 19, 2009;  Published online August 01, 2009}

\Abstract{Einstein gravity in both 3 and 4 dimensions, as well as some interesting genera\-lizations, can be written as gauge theories in which the connection is a Cartan connection for geometry modeled on a symmetric space. The relevant models in 3 dimensions include Einstein gravity in Chern--Simons form, as well as a new formulation of topologically massive gravity, with arbitrary cosmological constant, as a single constrained Chern--Simons action.  In 4 dimensions the main model of interest is MacDowell--Mansouri gravity, generalized to include the Immirzi parameter in a natural way.  I formulate these theories in Cartan geometric language, emphasizing also the role  played by the symmetric space structure of the model.  I also explain how, from the perspective of these Cartan-geometric formulations, both the topological mass in 3d and the Immirzi parameter in 4d are the result of non-simplicity of the Lorentz Lie algebra $\so(3,1)$ and its relatives.  Finally, I suggest how the language of Cartan geometry provides a guiding principle for elegantly reformulating any `gauge theory of geometry'.}


\Classification{22E70;
51P05;
53C80;
83C80;
 83C99}

\section{Introduction}

It is hard to overstate the importance of \'Elie Cartan's mathematical legacy, not only to mathe\-matics per se, but also to physics.  Cartan was surely among the earliest to contemplate the geometry of Einstein's spacetime in a deep way, and in this his superior geometric intuition was a distinct advantage.

Cartan was perhaps the f\/irst to appreciate that it is not the metric of spacetime that is most fundamental, but the connection (and the coframe f\/ield, which he introduced to the subject).  One piece of evidence for this is that in the Newtonian limit, while there is no sensible metric description of spacetime, Cartan found a beautiful reformulation of Newtonian gravity in the language of spacetime curvature, that is, in the language of connections \cite[Chapter~12]{MTW}.  Motivated by physical considerations, Cartan also generalized Riemannian geometry to include torsion \cite{Cartan-torsion}, and used this in a proposed generalization of Einstein gravity \cite{Trautman}.  This theory turned out to be important for coupling spinning particles to the gravitational f\/ield; though at the time spin had not been discovered in physics, spinors themselves were known~-- Cartan himself had discovered them in 1913~\cite{Cartan-spinors}.

Cartan's spaces with torsion include spaces with an absolute notion of parallelism between tangent spaces.  When Einstein proposed a unif\/ied theory based on absolute parallelism, a long and fascinating correspondence (1929--1932) between Cartan and Einstein ensued, concerning gravity, dif\/ferential equations, and geometry~\cite{EC}.

But, besides his work directly in relativity theory, Cartan seems to have had a knack for inventing exactly the sort of mathematics that turns out to be useful for spacetime physics.  Considering, for example, the Chern--Simons formulation of 3d general relativity \cite{AT,Witten3dgrav}, reviewed in Section~\ref{sec:einstein}, it is remarkable how many of the ingredients are due to Cartan.  The basic f\/ields are the coframe f\/ield $e$ and the spin connection~$\om$.  Cartan developed frame f\/ield methods extensively, and was probably the f\/irst to introduce them in general relativity; he was also the originator of connections in the general sense.
Moreover, the trick~-- rediscovered in physics many years later by MacDowell and Mansouri \cite{MM}~-- of combining the spin connection and coframe f\/ield into a single connection $A = \om \oplus e$, is an idea of Cartan's, at the foundation of the more general theory of Cartan connections.  A \cs\ action for this Cartan connection reduces to the Palatini action of 3d gravity:
\[
 S_{\scriptscriptstyle \rm CS}(A)
      = \int \tr\left( e\wedge {\star R} + \tfrac{1}{6}e\wedge {\star [e, e]}\right) .
\]
Here we notice more of Cartan's work: the action is built using Cartan's exterior dif\/ferential calculus~-- essential in all background-independent gravity theories, since, in the absence of any a priori metric structure, one needs a calculus that depends only on the dif\/ferentiable structure.  Finally, Cartan also invented the Killing form, written here as  `tr'.

In this expository article, I focus on just two aspects of Cartan's geometric work: symmetric spaces \cite{Helgason, KobayashiNomizu-2} and Cartan geometry \cite{AlekMic, Kobayashi-CC, Ruh, Sharpe}.   After brief\/ly reviewing the geometric ideas, I explain their application to certain gravity theories:
Einstein gravity in 3d and its extension to `topologically massive gravity', and 4d \MM\ gravity, generalized to include an Immirzi parameter.

\section{Reductive and symmetric Klein geometries}

An essential insight in Felix Klein's {\it Erlanger Programm} of 1872
was that problems in homogeneous geometry are `dual' to problems in group theory.  In particular, homogeneous manifolds are all isomorphic to $G/H$ for some Lie group $G$ and closed subgroup $H$, while `features' of a~given geometry can be understood by considering the subgroups of $G$ f\/ixing these features. In Klein's duality between geometry and algebra, `sub-features' and subgroups are inversely related in a Galois-like correspondence.

A standard example is Euclidean space $\R^n$.  It is homogeneous, meaning that it has a sym\-met\-ry group $G=\ISO(n):= \SO(n) \ltimes \R^n$ that acts transitively.  The subgroup stabilizing a~chosen point is $H\cong\SO(n)$, and so $\R^n$ may be identif\/ied with the coset space $\ISO(n)/\SO(n)$.  The subgroup stabilizing a line and a point on that line is isomorphic to $\O(n-1)$.  So, a typical geometric statement, like ``this point is contained in that line'',  has the algebraic counterpart ``this $\SO(n)$ subgroup contains that $\O(n-1)$ subgroup''.

Cartan was involved in both specialization and generalization of the Erlanger Programm, and these are the subject of this and the next section. Interestingly, both of these directions in Cartan's geometric research play a vital role in the formulations of gravity we shall discuss.  On one hand, he made signif\/icant contributions to the theory of `symmetric spaces', which have even `more' symmetry, in a certain precise sense, than generic Klein geometries.  On the other hand, he also laid the foundations for a radical generalization of Klein's work, in which the symmetries of an arbitrary smooth homogenous geometry hold only at an `inf\/initesimal' level, in much the same way that symmetries of Euclidean space only hold inf\/initesimally in Riemannian geometry.

In this section, I review some geometric and algebraic properties of special classes of Klein's homogeneous spaces:  {\em reductive spaces} and {\em symmetric spaces}.
Throughout this paper $G$ will denote a Lie group, $H\subseteq G$ a closed subgroup, which we think of as the stabilizer of a point in the homogeneous space $G/H$; the Lie algebras will be denoted $\g$ and $\h$.

Intuitively, a reductive space is one in which it makes sense, at least at an inf\/initesimal level, to ask whether a symmetry acts as a `pure translation' at a given point.  Let us explain this in a bit more detail. In an {\em arbitrary} homogeneous space, symmetries can be dif\/ferentiated, so that elements of the Lie algebra $\g$ (`inf\/initesimal symmetries') give tangent vectors at the basepoint of $G/H$ (`inf\/initesimal displacements').  A reductive space, roughly speaking, is one in which this process has a natural one-sided inverse; while there may be many inf\/initesimal symmetries giving rise to a given inf\/initesimal displacement, in a reductive geometry there is always a {\em canonical} choice of inf\/initesimal symmetry, the `pure translation' or `transvection'.  These transvections identify the tangent space at the basepoint of $G/H$ with a linear subspace~$\p$ of~$\g$.  To be consistent with symmetries, $\p$ must be invariant under the action of the stabilizer~$H$; this requirement automatically gives consistent behavior under a change of basepoint.

More formally, f\/irst observe that $\g$ is a representation of $\Ad(H)$ and has $\h$ as an invariant subspace.  When $\h$ has an $H$-invariant complement $\p$ in $\g$, the direct sum of $\Ad(H)$ representations
\begin{gather*}
\g= \h \oplus \p
\end{gather*}
is called a {\bf reductive splitting} of $\g$, and elements of  $\p$ are called {\bf transvections}.  The Klein geometry $G/H$ is {\bf reductive} if $\g$ admits such an invariant splitting.  Note in particular that for a reductive splitting we have
\[
      [\h,\h] \subseteq \h, \qquad [\h,\p] \subseteq \p.
\]

A symmetric space, intuitively speaking, is a space that looks the same if you {\em invert} it through any f\/ixed point.  Essentially, for an `internal observer', this means the space looks the same in any direction as it does in the opposite direction.

Formally, symmetric spaces may be def\/ined in a number of related ways.  The most general algebraic characterization, due to Cartan, is that a symmetric space is a homogenous space $G/H$ such that $H$ is the connected component of the set of elements of $G$ f\/ixed by an involutory group automorphism \cite{Chern-Chevalley}.  At the Lie algebra level, this means $\g$ has an involutory automorphism
\beq
\label{involution}
X \mapsto \widetilde X
\eeq
under which $\h$ is f\/ixed and has a complement $\p$ that is an eigenspace with eigenvalue $-1$, i.e.,
\[
    X = X_\h+ X_\p \quad \implies \quad \widetilde X = X_\h- X_\p.
\]
It is easy to check that this property is equivalent to the statement that $\g = \h \oplus \p$ is a~reductive splitting which also happens to be a $\Z/2$-grading of $\g$, where the subrepresentations $\h$ and $\p$ have respective grades 0 and 1, i.e.,
\[
      [\h,\h] \subseteq \h, \qquad [\h,\p] \subseteq \p, \qquad [\p,\p] \subseteq \h.
\]
A Lie algebra $\g$ equipped with such a $\Z/2$-grading is called a {\bf symmetric Lie algebra}.  (Of course, the $\Z/2$-grading here has no immediate relationship to supersymmetry.)  For purposes of this paper, we may def\/ine a {\bf symmetric space} to be a reductive geometry $G/H$ such that $\g= \h\oplus \p$ is a symmetric Lie algebra; this is the only property of symmetric spaces we shall need.    More information on symmetric spaces can be found in the references \cite{Helgason, KobayashiNomizu-2}.

Some of the most important homogeneous spaces used in physics are symmetric spaces.  In particular, the standard homogeneous solutions of Einstein's equations~-- de Sitter, Minkowski, and anti de Sitter spacetimes~-- are symmetric spaces; we describe these in detail in the next section, after a brief overview of Cartan geometry.

\section{Cartan geometry and Cartan gauge theory}
\label{sec:cartan-geom}

Cartan geometry \cite{Cartan1,Cartan2,Cartan3} has been called ``the most important among [Cartan's] works on dif\/ferential geometry'' \cite{Chern-Chevalley}.  It is also the foundation of some of the most elegant gauge-theoretic formulations of gravity, as I show in this paper.

In Cartan geometry, Klein geometry is generalized, roughly speaking, by considering spaces that look like $G/H$ only {\em infinitesimally}.   I sketch here what this means at a local level.  More details in the global setting, especially in the context of general relativity, may be found in~\cite{Wise}.  More complete mathematical references are the articles by Alekseevsky and Michor~\cite{AlekMic}, Kobayashi~\cite{Kobayashi-CC}, Ruh \cite{Ruh}, and the book by Sharpe \cite{Sharpe}.

Locally, a $G$-connection on $M$ is just a $\g$-valued 1-form
$A\maps TM \to \g$. A {\bf Cartan connection} satisf\/ies the additional requirement that the composite of $A$ with the quotient map:
\beq
\label{coframe}
\xymatrix{
    TM \ar[r]^{ A} \ar@/_3ex/[rr]_{e} & \g \ar[r] & \g/\h
}
\eeq
restricts to a linear isomorphism on each tangent space $T_xM$.  This obviously implies
\[
    \dim(M) = \dim(G/H).
\]
In fact, the $\g/\h$-valued 1-form $e$ in (\ref{coframe}), called the {\bf coframe f\/ield} or {\bf soldering form}, is viewed at each point as attaching, or `soldering', the tangent plane of $M$ to the tangent plane of the homogeneous model $G/H$.
Gauge transformations living in $H$ change the identif\/ication of $T_xM$ with $\g/\h$ via the adjoint representation, so any structure on $\g/\h$ that is $H$-invariant can be pulled back to the tangent spaces of $M$.  In particular, an invariant inner product on $\g/\h$ pulls back to give a semi-Riemannian metric on $M$.

The Cartan connection $A\maps TM \to \g$  itself can interpreted geometrically as giving instructions for `rolling' the homogeneous space $G/H$ along a path in $M$.  It associates, to each tiny motion along $M$, a tiny transformation of the homogeneous space that is `soldered' to it.

This idea is perhaps most vividly imagined by considering the example of Cartan geometry based on 2-dimensional spherical geometry, where $G= \SO(3)$, $H= \SO(2)$.  An instructive way to think of this example is by considering the motion of a `hamster-ball'~-- a transparent plastic ball whose motion over a surface $M$ is controlled by a hamster moving inside of it:
$$
\xy
(0,0)*{\includegraphics{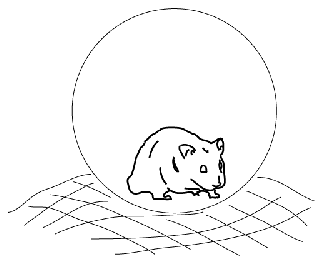}};
(15,8)*{G/H};
(18,-7)*{M};
\endxy
$$
Each tiny motion of the hamster in some direction tangent to the surface {\em completely determines} a~tiny rotation of the ball, an element of $\so(3)$.  This gives a Cartan connection locally described by a 1-form $A\maps TM\to \so(3)$.  If the hamster moves along some path $\gamma$ on the surface $M$, the holonomy of this Cartan connection along $\gamma$ then gives a group element describing how the ball rotates as it moves~\cite{Wise}.

The geometric beauty of Cartan's approach to connections is its concreteness.  While the more general principal bundle connections of Ehresmann are appropriate for the gauge theories of particle physics, whose gauge symmetries act on some abstract `internal space', Cartan connections are much more directly related to the geometry on the base manifold.

An example to keep in mind for Cartan geometry is the model Klein geometry $G/H$ itself.  The homogeneous space $G/H$ has a canonical Cartan connection, given by the Maurer--Cartan form $A\maps TG\to \g$, which is a 1-form on the total space of the principal $H$-bundle $G\to G/H$.  In local coordinates $A$ can be viewed as a $\g$-valued 1-form on $G/H$, and this local connection 1-form satisf\/ies the property that the coframe (\ref{coframe}) is nondegenerate.  If $\g/\h$ has an $H$-invariant inner product of some signature (for example, if the Killing form on $\g$ induces a nondegenerate form on the quotient), then pulling back along $e$ makes $G/H$ a homogeneous semi-Riemannian manifold.  General Cartan geometries modeled on $G/H$ may locally be viewed as deformations of $G/H$ itself.

If the model Klein geometry is {\em reductive}, then the quotient representation $\g/\h$ may be identif\/ied with a subrepresention $\p\subset \g$ of $H$, with $\g = \h\oplus \p$.  Given this splitting of $\g$, we can decompose $A$ uniquely as the a sum
\[
   A= \omega + e,   \qquad
        \omega \in \Omega^1(P,\h),\qquad
        e \in \Omega^1(P,\p)  ,
\]
where $\om$ is an $H$-connection and the coframe f\/ield is now a $\p$-valued 1-form.

The curvature of this connection is:
\beq
\label{reductive-curvature}
     F = R+\smallhalf[e,e] + d_\om e.
\eeq
The f\/irst term of this 2-form, the curvature $R= d\om+ \half[\om,\om]$ of the $H$-connection, takes values in $\h$, while the last takes values in $\p$, but the middle term, in a general reductive geometry, may have both $\h$ and $\p$ parts.  Thus
\[
    F_\h = R + \smallhalf [e,e]_\h \qquad F_\p =  d_\om e  + \smallhalf [e,e]_\p
\]
are the $\h$ and $\p$ components of the curvature,  called the {\bf corrected curvature} of  $\om$, and the {\bf torsion}.

Reductive geometries have many other nice properties.  For example, they have natural covariant derivative operators analogous to those of Riemannian geometry \cite[p.~197]{Sharpe}.   Also, the Bianchi identity in any reductive Cartan geometry splits neatly into $\h$ and $\p$ parts:
\[
     d_AF = 0  \qquad \iff \qquad d_\om R = 0  \qquad \text{and} \qquad d^2_\om e + [e,R] =0.
\]

If the reductive geometry happens to be {\em symmetric}, a number of additional special features emerge.  First, for a Cartan connection $A= \om +e$ based on a symmetric geometry, $[e,e]$ is a~2-form with values in $\h$, so the curvature (\ref{reductive-curvature}) splits more cleanly:
\[
    F_\h = R + \smallhalf [e,e],  \qquad F_\p = d_\om e.
\]
Moreover, since $\g = \h \oplus \p$ is a symmetric Lie algebra, the connection $A$ can be composed with the involution \eref{involution} to give a new connection:
\[
        \widetilde A = \om - e.
\]
The curvature of $\widetilde A$ is:
\[
       \widetilde F = R+\smallhalf[e,e] - d_\om e,
\]
so the map $A\mapsto F[A]$ taking a connection to its curvature is equivariant with respect to the involution (\ref{involution}) on $\g$.
This give a nice way of calculating  $\h$ and $\p$-valued parts of the curvature in the symmetric case:  the part of $F$ symmetric under the involution:
\[
F_\h = \smallhalf(F+\widetilde F) = R+\smallhalf[e,e]
\]
is the corrected curvature, while the antisymmetric part:
\[
F_\p  = \smallhalf(F-\widetilde F) = d_\om e
\]
is the torsion.

\subsection*{de Sitter, Minkowski, and anti de Sitter models}

Let us now give more explicit descriptions of Cartan geometries for the most important symmetric space models for studying gravity, namely the homogeneous constant-curvature spacetimes with cosmological constant $\Lambda$.  For spacetime dimension $n$, the relevant symmetry groups are:
\beq
\label{sign}
    G = \left\{
      \begin{array}{lll}
          \SO(n,1), & \rule{1em}{0pt} \Lambda > 0  &  \text { (de Sitter)},\\
          \ISO(n-1,1),& \rule{1em}{0pt} \Lambda = 0 & \text { (Minkowski)}, \\
          \SO(n-1,2), & \rule{1em}{0pt} \Lambda < 0 &\text { (anti de Sitter)}
      \end{array}
      \right.
\eeq
and in each case the stabilizer is the Lorentz group $H = \SO(n-1,1)$.  These are each symmetric Lie algebras:
\beq
\label{split}
     \g = \h \oplus \p,
\eeq
where $H$ acts as Lorentz transformations on $\p\iso \R^n$.

A Cartan connection $A$ modeled on one of the homogeneous spacetimes $G/H$ above may be written in index notation as ${A^I}_J$, where the indices $I,J = 0,1,\ldots n$ label matrix components in the fundamental representation of $\g$.   Letting lower case indices $a, b, \ldots$, take values $0, 1, \ldots, n-1$, we can take the $\h$ part of $A$, a Lorentz connection, to be the upper left block of matrix components of $A$:
\[
         {A^a}_b={\om^a}_b.
\]
This leaves the last row and column for the components of the $\p$ part, which we write
\[
       {A^a}_n =  \frac{1}{\ell}e^a  ,   \qquad   {A^n}_a =  - \frac{\xe}{\ell}e_a.
\]
Here we scale the coframe f\/ield $e$ by a constant $\ell$ with units of length, and $\xe$ is the sign of the cosmological constant:
\[
    \xe = \sgn(\Lambda).
\]
 We can then calculate the $\h$ and $\p$ components of the curvature $F= dA + A\wedge A$:
\[
     {F^a}_b     = {R^a}_b  - \frac{\xe}{\ell^2}\, e^a \wedge e_b
\qquad \quad
    {F^a}_{n}   = \frac{1}{\ell}\; d_\om e^a.
\]
The f\/irst of these components explains the term `corrected curvature' for $F_\h$: if we constrain the torsion part $d_\om e$ to vanish and choose, in dimension $n \ge 3$,
\[
     \frac{\xe}{\ell^2} = \frac{2\Lambda}{(n-1)(n-2)},
\]
then the condition for the corrected curvature to vanish is precisely that spacetime be locally isometric to the homogeneous model with cosmological constant $\Lambda$.  Indeed, rewritten in metric variables, the corrected curvature in any of these examples is
\[
  {\mathcal R}_{\mu\nu\rho\sigma}=    R_{\mu\nu\rho\sigma} - \frac{2\Lambda}{(n-1)(n-2)} \left(g_{\mu\rho}g_{\nu\sigma} - g_{\mu\sigma}g_{\nu\rho}\right).
\]
The correction term is precisely the Riemann tensor of the model $G/H$, so $R_{\mu\nu\rho\sigma}$ satisf\/ies the cosmological Einstein equation ${G}_{\mu\nu} + \Lambda g_{\mu\nu} = 0$ if and only if ${\mathcal R}_{\mu\nu\rho\sigma}$ satisf\/ies the `cosmologically corrected' Einstein equation ${\mathcal G}_{\mu\nu} = 0$, where the cosmological Einstein tensor $\G_{\mu\nu} =  {\mathcal R}_{\mu\nu} - \half { {\mathcal R} g_{\mu\nu}}$ is built from $ {\mathcal R}_{\mu\nu} = { {\mathcal R}^\rho}_{\mu\rho\nu}$ and $ {\mathcal R} =  {\mathcal R}^\mu_\mu$, the corrected Ricci tensor and scalar curvature.

Translating between the Lie algebra-valued form language and the index notation requires a bit of care, especially with with respect to signs.  In particular, the expression $\half[e,e]$ for the $\p$-valued 1-form $e$ means the same as $e^a \wedge e^b$ or $-e^a \wedge e^b$ or $0$, depending on the sign of $\Lambda$, as indicated above.  However, the form notation is tremendously unifying. In the rest of this paper, I use only the Lie algebra-valued dif\/ferential form notation, and use units
\beq
\label{unit}
  \ell=1,
\eeq
writing simply $A = \om + e$, as in the general discussion above.  This has the advantage that, in the torsion-free case, the equation
\[
     R+ \tfrac 12 [e,e] = 0
\]
describes locally homogeneous geometry for {\em any} of the spacetime models we have just described; all information about signs is hidden in the Lie bracket, according to the choice of gauge group~(\ref{sign}), and for many purposes we may thus treat all cases in parallel.

Elsewhere \cite{Wise} I have emphasized that MacDowell--Mansouri gravity is a gauge theory for a~Cartan connection based on these models.  In fact, there are many more gravity models that are best viewed in this way.  I review a few of these in the following sections.  (See also the references \cite{Aldinger, HMMN, MOH, Zardecki} for related applications of Cartan geometry to gravity.)

\section{3d Einstein gravity}
\label{sec:einstein}

The f\/irst clue that 3d gravity is a gauge theory for a Cartan connection is quite simple.  In the language of the previous section, the 3d Einstein equations may be written as
\[
      R(\om)+\tfrac 12 [e,e] = 0, \qquad d_\om e = 0,
\]
where $\om$ and $e$ are the components of a Cartan connection $A=\om + e$ for the  gauge group (\ref{sign}) of gravity (with $n=3$) appropriate to the cosmological constant, and $R$ is the curvature of $\om$.   These are summarized by the single equation
\[
       F = 0.
\]
which simply says the Cartan connection is `f\/lat', meaning the Cartan geometry is locally the same as that of the Klein model.

In fact, solutions of 3d gravity have been studied using a special case of Cartan geometry, though not in quite this language.   For f\/lat Cartan geometries, the `rolling' described in Section~\ref{sec:cartan-geom} is essentially trivial and spacetime is locally isometric to the Kleinian model spacetime $G/H$.  Flat Cartan geometry reduces precisely to the theory of {\bf geometric structures}.  A manifold with a $G/H$ geometric structure can be viewed as a manifold whose charts are maps into the homogeneous model $G/H$, rather than $\R^n$, and whose transition functions are required to live in the symmetry group $G$.
These were descibed by Ehresmann \cite{Ehresmann2}, and studied extensively by Thurston \cite{Thurston, Thurston2}, particularly in the 3d Riemannian case, where they are the basis for the geometrization conjecture, now proven, largely by Perelman.  The theory of geometric structures is an important tool in geometric toplogy (see,  e.g.~Goldman's papers \cite{Goldman2,Goldman3} and references therein).  It has been used by Carlip to study solutions of 3d general relativity, and 3d quantum gravity \cite{Carlip, Carlip2}.

This Cartan-geometric formulation of 3d general relativity is ultimately based on a \cs\ theory for a Cartan connection.  Gravity in 3 dimensions was studied as a \cs\ theory f\/irst by Ach\'ucarro and Townsend \cite{AT}, and then in further detail by Witten \cite{Witten3dgrav}.  While these ideas are in principle familiar, here I cast the Chern--Simons formulation in the language of Cartan connections for the relevant symmetric Lie algebras, and indicate what role is played by the symmetric space structure of the model.
In particular, in Proposition~\ref{2CS}, we shall see that the \cs\ formulation can be elegantly expressed in terms of the symmetric space structure of the de Sitter, Minkowski, and anti de Sitter spacetimes:  3d general relativity is the part of \cs\ theory that is antisymmetric with repect to the involution on the model symmetric space.

Given any nondegenerate invariant bilinear form $\bilin\maps \g\tensor \g \to \R$, the Chern--Simons action
 \[
  S_{\rm \scriptscriptstyle CS}^\bilin(A) = \tfrac 12 \int \bilin (A\wedge dA + \tfrac{1}{3} A\wedge [A,A])
 \]
results in the classical f\/ield equations $F=0$, which may be interpreted as the 3d Einstein equations.  But what the \cs\ action for $A$ looks like in terms of the constituent f\/ields $\om$ and $e$ depends on the properties of the bilinear form $\bilin$, as the following straightforward propositions show.
\begin{proposition}
\label{CS-null}
Let $G/H$ be a symmetric space such that the symmetric splitting $\g= \h \oplus \p$ has $\h\perp \h$ and $\p\perp \p$ for an invariant bilinear form $\bilin\maps \g\tensor \g \to \R$.
If $A = \om + e$ is a symmetric Cartan connection, then up to boundary terms,
\[
   S_{\rm \scriptscriptstyle CS}^\bilin(A) =
   \int \bilin\left( e\wedge R + \tfrac{1}{6}e\wedge [e, e]\right).
\]
\end{proposition}

\begin{proposition}
\label{CS-perp}
Let $G/H$ be a symmetric space such that the symmetric splitting $\g= \h \oplus \p$ has $\h\perp \p$ for an invariant bilinear form $\bilin\maps \g\tensor \g \to \R$.  If $A = \om + e$ is a symmetric Cartan connection, then up to boundary terms,
 \begin{gather*}
S_{\rm \scriptscriptstyle CS}^\bilin(A)
    =   \int \bilin \left(\om\wedge \left(d\om + \tfrac{1}{3} [\om,\om]\right) + e\wedge d_\om e\right).
 \end{gather*}
\end{proposition}

If $G$ is a simple group, then every invariant inner product on $\g$ is a multiple of the Killing form, so in this case there is only one choice of $\bilin$ up to normalization.  The choice of normalization corresponds to a choice of coupling constant, which inf\/luences the strength of possible interaction terms and scales quantum f\/luctuations.

The situation in 3d gravity is quite dif\/ferent: the relevant gauge groups:
\[
G = \left\{
      \begin{array}{ll}
          \SO(3,1), & \rule{1em}{0pt} \Lambda > 0,  \\
          \ISO(2,1),& \rule{1em}{0pt} \Lambda = 0, \\
          \SO(2,2), & \rule{1em}{0pt} \Lambda < 0.
      \end{array}
      \right.
\]
 are not simple.  In fact, the Lie algebras each have {\em two}-dimensional spaces of invariant inner products.  In particular, as observed by Witten, for $\so(3,1)$ and $\so(2,2)$, it is possible to choose the inner product $\bilin$ such that the conditions of either Proposition~\ref{CS-null} or Proposition~$\ref{CS-perp}$ are satisf\/ied. The obvious choice, obtained by taking $\bilin$ proportional to the Killing form, corresponds to the situation of Proposition~\ref{CS-perp}.  But, we may also insert an `internal' Hodge star into one of the factors~-- the Lie algebra $\g$ for $\Lambda\neq 0$ has a notion of Hodge duality by virtue of the isomorphism $\g\cong\Lambda^2\R^4$~-- and this puts us in the situation of Proposition~\ref{CS-null}.  In fact, the inner product $\tr(X{\star Y})$ makes sense for $\Iso(2,1)$, as well, and is nondegenerate, and it is this inner product that gives the 3d Einstein--Palatini action with $\Lambda =0$~\cite{Witten3dgrav}.

Explicitly, letting `tr' denote the Killing form (or a multiple thereof), these two choices lead to the distinct actions:
\begin{gather*}
S_{\rm \scriptscriptstyle CS}^{\tr\star}(A) : = \tfrac 12 \int \tr\left(A\wedge {\star\!}\left(dA + \tfrac{1}{3} [A,A]\right)\right) \qquad \text{and} \\
  S_{\rm \scriptscriptstyle CS}^{\tr}(A) : =   \tfrac 12 \int \tr\left(A\wedge \left(dA + \tfrac{1}{3} [A,A]\right)\right),
\end{gather*}
where the latter is nondegenerate for $\Lambda\neq 0$.
The version with the Hodge star gives the Palatini action for 3d general relativity:
\beq
\label{CS-star}
S_{\rm \scriptscriptstyle CS}^{\tr\star}(A)  =
   \int \tr\left( e\wedge {\star R} + \tfrac{1}{6}e\wedge {\star [e, e]}\right)  ,
\eeq
while omitting the star gives \cs\ in $\om$ plus the torsion term:
\beq
\label{CS-nostar}
S_{\rm \scriptscriptstyle CS}^{\tr}(A) =
   \int \tr\left(\om\wedge \left(d\om + \tfrac{1}{3} [\om,\om]\right) + e\wedge d_\om e\right).
\eeq

Varying $e$ and $\om$ independently, the resulting variations of these actions, up to boundary terms, are respectively
\[
   \int \tr\left(\xd \om \wedge (R+\tfrac 12 [e,e]) + \xd e \wedge d_\om e\right)
\qquad
\text{and}
\qquad
   \int \tr\left(\xd e \wedge \star (R + \tfrac 12 [e,e]) + \xd \om \wedge {\star d_\om e}\right).
\]
The equations of motion resulting from either of these actions are thus the same, as expected, since both must ultimately say ``$F=0$''.  However, it is worth observing that in (\ref{CS-star}) the Einstein equation comes from variation of $e$ and the vanishing of the torsion from the variation of $\om$, while in (\ref{CS-nostar}) these roles are reversed.  While the classical theories corresponding to (\ref{CS-star}) and (\ref{CS-nostar}) are the same, they dif\/fer as quantum theories~\cite{Witten3dgrav}.

More generally, we get a theory classically equivalent to Palatini-style 3d general relativity by taking a generic linear combination of (\ref{CS-star}) and (\ref{CS-nostar}).  This corresponds to using a Chern--Simons action with a generic bilinear form
\beq
\label{innerprod}
    \bilin(X,Y) = \tr(X(c_0 + c_1{\star})Y),
\eeq
where $c_0,c_1\in \R$ are constants.  Every invariant bilinear form on the Lie algebra $\g$ is of this form, as shown in Appendix \ref{apx:Lie-alg}.
\begin{proposition}
\label{Einstein-CS}
Let $\g = \so(3,1)$, $\Iso(2,1)$, or $\so(2,2)$.  Starting with any invariant symmetric bilinear form $\bilin$ on $\g$, and defining the \cs\ action
\[
  S_{\rm \scriptscriptstyle CS}^{\bilin}(A) = \tfrac 12 \int \bilin (A\wedge dA + \tfrac{1}{3} A\wedge [A,A])
\]
we have
\[
  S_{\rm \scriptscriptstyle CS}^{\bilin}(A) = c_1 S_{\rm Pal}(\om,e) + c_0 S_{\rm \scriptscriptstyle CS}(\om)  + c_0 \int\tr(e\wedge d_\om e),
\]
where $(c_0,c_1)$ are constants that relate the inner product $\bilin$ to the Killing form $\tr$ via \eqref{innerprod}.
\end{proposition}

Since $\tr(XY)$ is nondegenerate in the $\Lambda\neq 0$ cases, $\bilin$ will be degenerate if and only if $(c_0 + c_1{\star})Y = 0$ for some $Y\in \g$.  This happens in the anti de Sitter case, where the choice $c_0=\pm c_1$ annihilates (anti) self dual elements of $\so(2,2)$.  For $\Lambda=0$, $\bilin$ is degenerate on $\Iso(2,1)$ if  $c_1 = 0$.

Proposition~\ref{Einstein-CS} shows we can pick the inner product $\bilin$ on $\g$ to obtain any linear combination of Palatini and \cs\ actions we wish.  A related question is whether, given a generic inner product $\bilin$, we can somehow recover the standard Palatini action.  This can be done if we make use of the involution (\ref{involution}) for our  symmetric Lie algebra.  Composing this involution with the connection $A$ maps $e\mapsto -e$, and the action (\ref{CS-star}) changes sign, while (\ref{CS-nostar}) is invariant.  This easily implies Palatini gravity is the ``antisymmetric part'', with respect to the involution on $\g$, of a generic \cs\  action.  More precisely, this shows:

\begin{proposition}
\label{2CS}
Let $\g = \so(3,1)$, $\Iso(2,1)$, or $\so(2,2)$.  Starting with any invariant symmetric bilinear form $\bilin$ on $\g$, and defining the \cs\ action
\[
  S_{\rm \scriptscriptstyle CS}^{\bilin}(A) = \tfrac 12 \int \bilin \left(A\wedge dA + \tfrac{1}{3} A\wedge [A,A]\right)
\]
we have
\[
\tfrac 12 S_{\rm \scriptscriptstyle CS}^\bilin(A)
+\tfrac 12 S_{\rm \scriptscriptstyle CS}^\bilin(\widetilde A)
= c_0 S_{\rm \scriptscriptstyle CS}(\om) + c_0 \int \tr(e\wedge d_\om e)
\]
and
\[
\tfrac 12 S_{\rm \scriptscriptstyle CS}^\bilin(A)
-\tfrac 12 S_{\rm \scriptscriptstyle CS}^\bilin(\widetilde A)
= c_1    \int \tr\left( e\wedge {\star R} + \tfrac{1}{6}e\wedge {\star [e, e]}\right),
\]
where $(c_0,c_1)$ are constants that relate the inner product $\bilin$ to the Killing form $\tr$ via \eqref{innerprod}.
\end{proposition}

It seems worth remarking here that Proposition~\ref{2CS} is at least superf\/icially dif\/ferent from the most familiar way of writing 3d gravity as a dif\/ference of \cs\ theories.  The standard approach \cite{Witten3dgrav}, in the anti de Sitter case $\Lambda<0$, is based on the Lie algebra isomorphism
\beq
\label{coincidence}
      \so(2,2) \iso \Sl(2,\R)\oplus \Sl(2,\R)
\eeq
obtained by splitting $X\in \so(2,2)$ into a self dual part $X^\ssp$ and an anti-self dual part $X^\ssm$:
\[
     X = X^\ssp + X^\ssm, \qquad X^\sspm = \smallhalf(X \pm {\star X}), \qquad {\star X^\pm} = \pm X^\sspm.
\]
The self dual and anti self dual
subalgebras are Killing-orthogonal, and rewriting Proposition~\ref{2CS} in terms of these gives the usual $\Sl (2,\R)\times \Sl(2,\R)$ version.
Unfortunately, this approach is less satisfactory for $\Lambda>0$, since $\so(3,1)$ has a splitting analogous to (\ref{coincidence}) only after complexif\/ication, forcing one to use a complex connection $A = \om + {\rm i}e$, the geometric meaning of which is unclear. The approach used in this section, based on a symmetric Cartan connection, is more natural: it puts the $\so(3,1)$ and $\so(2,2)$ versions on equal footing, since what is important here is not the Lie algebra coincidence (\ref{coincidence}), but the symmetric splitting (\ref{split}) that all of the relevant gravity gauge groups have.

\section{Topologically massive gravity}

It is indeed surprising that the classical equations of 3d general relativity arise from any member of an inf\/inite family of inequivalent actions, corresponding to dif\/ferent choices of invariant inner product (\ref{innerprod}) on $\g$ (Witten \cite{Witten3dgrav} aptly described this behaviour as ``enigmatic").  We can shed light on this situation by considering the extension of 3d general relativity to `topologically massive gravity', where we shall see that the classical theory {\em does} depend on the choice of inner product in the \cs\ action.  In fact, here we f\/ind that the 2-dimensional space of invariant inner products correspond to the two coupling parameters in the theory.

In topologically massive gravity \cite{Deser:1982vy}, one takes (the negative of) the Einstein--Hilbert action and adds a Chern--Simons term in the spin connection:
\beq
\label{tmg-action}
\displaystyle
S_{\scriptscriptstyle\rm TMG}(e) = \int   \tr \left(
- e\wedge {\star R} - \tfrac{1}{6}e\wedge{\star [e,e]}
+
\tfrac{1}{2\mu}
 \left(\om\wedge
d\om
+\tfrac{1}{3}
\om \wedge
[\om, \om] \right) \right),
\eeq
where the coupling parameter $\mu$ is called the `topological mass'.   We emphasize that this action depends only on the coframe f\/ield $e$, and $\om = \om[e]$ is the torsion-free spin connection.
\beq
\label{tmg-eom}
G + {\Lambda} e + \tfrac 1{\mu}
     C = 0,
\eeq
where $G$ is the Einstein tensor, $C$ is the Cotton tensor~-- the conformal curvature in three dimensions~-- and the cosmological constant $\Lambda$ is $+1$, $0$ or $-1$, according to our choice of units~(\ref{unit}).

Paradoxically, while standard 3d Einstein gravity is a topological f\/ield theory, adding the `topological mass' term in (\ref{tmg-action}) makes the theory {\em not} topological~-- solutions of the equations of motion (\ref{tmg-eom}) are not all locally gauge equivalent.  The somewhat misleading name `topologically massive' comes from the topological character of Chern--Simons theory and because adding this term here gives {\em mass} to the gauge f\/ield.  In fact, the linearized approximation of the theory may be described in terms of massive scalar wave equations \cite{CDWW,Deser:1982vy}.

It follows immediately from Proposition~\ref{Einstein-CS} that topologically massive gravity (\ref{tmg-action}) is simply a~single \cs\ action, once we impose the zero-torsion constraint $d_\om e = 0$:
\begin{proposition}
Let $\g = \so(3,1)$, $\Iso(2,1)$, or $\so(2,2)$, and form a Cartan connection from just the coframe field, defined by $A(e) = \om(e) + e$, where $\om(e)$ is the torsion-free spin connection.
Then the action for topologically massive gravity, with cosmological constant appropriate to $\g$, is given by
\beq
\label{CS-tmg}
S_{\scriptscriptstyle\rm TMG}(e) =
 S_{\scriptscriptstyle\rm CS}^\bilin(A(e)),
\eeq
where $S_{\scriptscriptstyle\rm CS}^\bilin$ is the \cs\ action for the bilinear form $\bilin(X,Y) = \tr(X(\frac{1}{\mu} + \star)Y)$.
\end{proposition}

When $\Lambda <0$, so that $\g= \so(2,2)$, the choice $\mu = \pm 1$ makes the \cs\ action in this proposition degenerate.   However, it becomes nondegenerate if we restrict to (anti) self dual connections $A= \mp {\star A}$.  Interestingly, these particular values $\mu = \pm 1$ of the mass have been the subject of much recent work in 3d gravity, sparked by a claim \cite{Li} that the theory, while nontrivial, becomes `topological in the bulk' at these values.  What precisely happens at these values of $\mu$ is still somewhat mysterious; I shall not consider this topic further here (see \cite{CDWW} for many references), but focus on the geometric content of the theory for generic $\mu$.

One can also ask, in parallel with Proposition~\ref{2CS}, how to construct \tmg\ starting with any generic invariant inner product on $\g$.  This is easy, using the results of Proposition~\ref{2CS} and imposing the $d_\om e = 0$ constraint: if we normalize our bilinear form so that $\bilin(X,Y) = (X(c_0 + {\star})Y)$, then the result is
\beq
\label{2CS-tmg}
S_{\scriptscriptstyle\rm TMG}(e) =
 -\tfrac 12 \left(1-{\tfrac 1{\mu c_0}}\right)S_{\scriptscriptstyle\rm CS}(A(e))
  + \tfrac 12 \left(1+{\tfrac 1{\mu c_0}}\right) S_{\scriptscriptstyle\rm CS}(\widetilde A(e)).
\eeq
This generalizes to arbitrary $\Lambda$ (and to arbitrary $\bilin$) our earlier result \cite{CDWW} that topologically massive gravity with $\Lambda < 0$ is a linear combination of Chern--Simons theories.  There we used the approach based on the Lie algebra coincidence~(\ref{coincidence}) and wrote topologically massive AdS$_3$ gravity as a sum of Chern--Simons actions whose coef\/f\/icients
\[
\left(1+\tfrac{1}{\mu}\right) \qquad \text{and} \qquad \left(1-\tfrac{1}{\mu}\right)
\]
turn out to coincide with the left and right central charges of the boundary conformal f\/ield theory.   Rewriting (\ref{2CS-tmg}) in the $\Sl(2,\R)\oplus \Sl(2,\R)$ basis recovers this result.

A few words about using a f\/irst-order formalism in topologically massive gravity are in order.  Deser \cite{Deser} has analyzed the obstacles to using the Palatini approach in the original topologically massive gravity action (\ref{tmg-action}).  Here it was pointed out that there are really four choices: one can choose, independently in the Einstein term and the \cs\ term, either to use an arbitrary connection $\om$, or the Levi-Civita spin connection $\widetilde \om(e)$.  The interesting results are:
\begin{enumerate}
  \item Forcing the torsion constraint only in the Einstein term, while using an arbitrary  connection in the \cs\ term:
  \[
    S(\om,e) :=  S_{\scriptscriptstyle \rm EH}(\widetilde \om(e)) + \tfrac{1}{\mu}S_{\scriptscriptstyle \rm CS}(\om)
  \]
recovers pure metric Einstein gravity.
  \item Allowing an arbitrary connection in the Einstein term (as in the original Palatini approach), while imposing the torsion constraint on the \cs\ term:
  \[
    S(\om,e) :=  S_{\scriptscriptstyle \rm EH}(\om,e) + \tfrac{1}{\mu}S_{\scriptscriptstyle \rm CS}(\widetilde\om(e))
  \]
recovers topologically massive gravity.
\end{enumerate}
Using the Levi-Civita connection in both terms, one has the standard action (\ref{tmg-action}) for topologically massive gravity, while the f\/inal possibility~-- independent coframe and connection in both terms~-- was found to be much less tractable~\cite{Deser}.  This is perhaps best explained by considering (\ref{CS-tmg}) as the def\/ining action of topologically massive gravity, rather than~(\ref{tmg-action}).  The two actions~(\ref{tmg-action}) and~(\ref{CS-tmg}) are equivalent when the torsion constraint is imposed.  However, if we relax this constraint and use a Palatini approach, they are not: they dif\/fer by the term
\[
       \int \tr (e \wedge d_\om e).
\]
From the perspective of (\ref{CS-tmg}), then, a naive application of the Palatini method to topologically massive gravity gives, not a f\/irst-order version of topologically massive gravity, but simply ordinary 3d Einstein gravity.

A genuine f\/irst order formulation of topologically massive gravity may be obtained by adding to the action a Lagrange multiplier term forcing the torsion~-- or in other words, the $\p$-part of the curvature, $F_\p = \half(F-\widetilde F)$~-- to vanish.  A similar f\/irst order formulation has been used to analyze the constraint algebra of topologically massive gravity \cite{Carlip-constraints}.  However, the latter approach uses an $\SL(2,\R)$ connection $\om + \mu{\star e}$, with $\mu$ as the length scale, as opposed to the $G$ connection $\om \oplus e$ used here, and as a result involves an additional constriant term.

\section[MacDowell-Mansouri gravity]{MacDowell--Mansouri gravity}

The action of MacDowell--Mansouri gravity \cite{MM,Wise} is:
\begin{equation}
\label{mm-action}
  S_{\scriptscriptstyle \rm MM}(A) = - \tfrac{1}{2} \int \tr(F_\h \wedge \star F_\h).
\end{equation}
Here $F$ is the curvature of a Cartan connection $A= \om + e$ for the 4d de Sitter or anti de Sitter model, so the relevant gauge groups are
\[
\label{sign}
H= SO(3,1),
\qquad
    G = \left\{
      \begin{array}{ll}
          \SO(4,1), & \rule{1em}{0pt} \Lambda >  0, \\
          \SO(3,2), & \rule{1em}{0pt} \Lambda < 0.
      \end{array}
      \right.
\]
While this looks vaguely like a topological theory of the form $\int F^2$, two crucial tricks make the theory nontrivial: the projection $F\mapsto F_\h$ into the $\so(3,1)$-part of $F$:
\[
   F_\h = \tfrac 12 (F + \widetilde F)
\]
and also the Hodge star operator $\star$ on $\so(3,1)$.     As shown below, the action reproduces general relativity.  Since I have given a detailed account of the Cartan-geometric underpinnings of this theory elsewhere~\cite{Wise}, I focus here on a couple of features not brought out there.

First, a few brief comments about gauge invariance are in order.  One initially unattractive feature of the MacDowell--Mansouri formulation, from a gauge theory perspective, is that in the action (\ref{mm-action}) we seem to have broken the $G$-invariance `by hand' down to some specif\/ic $\SO(3,1)$ subgroup.  This was certainly the implication in the original \MM\ paper, where an index in the action was explicitly f\/ixed to a specif\/ic value.  However, one can also set up \MM\ gravity so that the de Sitter symmetry is spontaneously broken.  The selection of an $\SO(3,1)$ subgroup is a gauge choice, in that it can be made independently at each point, and the choice has no bearing on the physics. This choice amounts to picking a~basepoint in the Klein geometry associated to each point \cite{Wise}.  The idea that this ``section of basepoints'' should serve as a gauge f\/ield in Cartan-type gauge theory is noting new~-- for example, twenty years ago, Aldinger \cite{Aldinger} proposed a related $\SO(4,1)$ Cartan-type gauge theory, in which this section serves as a Higgs f\/ield.  For more recent discussions of spontaneous symmetry breaking in MacDowell--Mansouri-like gauge theories of gravity, see e.g.\ \cite{FreidelStarodubtsev, GibbonsGielen, Ikeda}.

Now, a key point is that, even though (\ref{mm-action}) is based on four-dimensional Cartan geometry, the action really only uses the Killing form on the stabilizer $\h =\so(3,1)$.  This is the same as the Lie algebra generating isometries of {\em three}-dimensional de Sitter space.  We can thus use the same trick as in 3d gravity, and try using a more general inner product
\[
     \bilin(X,Y) = \tr(X(c_0 + c_1{\star})Y)
\]
on the Lie algebra $\h$, giving the generalized \MM\ action
\begin{equation*}
  S^\bilin_{\rm\scriptscriptstyle  MM}(A) = - \tfrac{1}{2} \int \bilin(F_\h \wedge F_\h).
\end{equation*}

Substituting $F_\h = R + \half [e,e]$ in this action leads to a variety of terms.  The terms involving
\[
   \tr(R\wedge R) \qquad \text{and}\qquad    \tr(R\wedge {\star R})
\]
are topological invariants with vanishing variation, due to the Bianchi identity and the identity $\xd R = d_\om \xd \om$. These are well-known to arise in the MacDowell--Mansouri approach, and have been discussed elsewhere (see e.g.~\cite{FreidelStarodubtsev, NSO}); here I simply discard them for convenience.  Also, the Killing form on $\so(3,1)$ extends to the Killing form on all of $\g$, and this lets us use adjoint-invariance to conclude that
\[
     \tr ([e,e]\wedge [e,e]) = \tr(e\wedge [e,[e,e]])
     = 0
\]
since $[e,[e,e]]$ vanishes by the graded Jacobi identity for $\g$-valued forms.  Fortunately, the Hodge star saves the other $e^{ 4}$ term from the same fate, giving the cosmological constant term.   We are thus left with:
\[
   S^\bilin_{\rm\scriptscriptstyle  MM}(A) =  - \int \tr\left([e,e]\wedge{\star R}+ \tfrac{1}{4} [e,e] \wedge \star[e,e]) + \tfrac{1}{\gamma} [e,e]\wedge R \right),
\]
where we have identif\/ied $(c_0,c_1) = (1,1/\gamma)$.  Noting that in 4 dimensions our choice of units $\ell = 1$ in (\ref{unit}) implies $\Lambda = \pm 3$, this is the action for general relativity with {\bf Immirzi parameter}~$\gamma$.  This may be compared to the action of Freidel and Starodubtsev~\cite{FreidelStarodubtsev}, which uses a BF theory action for $G$, `perturbed' by a $B\wedge B$ term which essentially uses this generalized inner product.  The generalization of MacDowell--Mansouri gravity to include the Immirzi parameter has also been nicely discussed by Randono~\cite{Randono}.

\section{Conclusion}

Starting from the 4d Palatini action and interpreting the spin connection and coframe f\/ield as components of a Cartan connection, one is led almost automatically to the \MM\ formulation.  (This was my viewpoint in~\cite{Wise}.)  A similar statement can be made about passing from the 3d Palatini action to the Chern--Simons formulation.  These examples suggest a general strategy for neatly rewriting geometric f\/ield theories in the geometric language of Cartan connections.  Since Cartan geometry is the natural connection-based language for describing quite general sorts of geometry, one might suspect that this should be possible for any sensible `gauge theory of geometry'.

As an application of this principle, I have here given a new formulation of topologically massive gravity as a {\em single} constrained Chern--Simons theory.  This formulation is related to the $\SL(2,\R)\times \SL(2,\R)$ formulation given previously \cite{CDWW}, though there a linear combination of two Chern--Simons theories was necessary.

Some properties used in these Cartan-geometric gravity theories are special to de Sitter, Minkowski, and anti de Sitter models.  In particular, the Hodge star exists only on certain special Lie algebras.  On the other hand, many calculations in this paper have used nothing particular to these models, but would work for Cartan geometry based on any symmetric Klein geometry $G/H$.  It would be interesting to generalize these theories, writing down actions for geometric gauge theories based on other Cartan geometries.

\appendix

\section{The modesty of Lie algebra-valued dif\/ferential forms}
\label{apx:forms}

In the spirit of \'Elie Cartan, I have avoided as much as possible ``the debauch of indices'', using instead the cleaner language of Lie algebra-valued forms.
Here I provide, especially for convenience of the reader more accustomed to index-based calculations, some of the standard formulas.
All basic calculations in the f\/ield theories described in this paper may be carried out using the tools presented here.

For $\g$ a Lie algebra, denote the space of $\g$-valued $p$-forms on $M$ by
$\Omega^p(M,\frakg) := \Omega^p(M)\tensor \frakg.
$
The bracket of $\g$-valued forms is def\/ined using the wedge product on form parts and the Lie bracket on Lie algebra parts:
\[
    \om = \om^\xa \tensor v_\xa, \;\mu = \mu^\xb \tensor v_\xb \quad \implies \quad
   [\omega,\mu] := (\omega^\xa \wedge \mu^\xb) \tensor [v_\xa, v_\xb],
\]
where $\{v_\xa\}$ is a basis of $\g$.
Like the wedge product of ordinary forms, this bracket is graded-commutative, but with an additional minus sign coming from the Lie bracket:
\begin{equation*}
[\omega,\mu] = (-1)^{pq+1}[\mu,\omega],
\qquad
\om \in \Omega^p(M,\g),\qquad  \mu \in \Omega^q(M,\g).
\end{equation*}
It also obeys a graded Jacobi identity:
\[
     [\lambda, [\om,\mu]] = [[\lambda,\om],\mu] + (-1)^{lp}[\om,[\lambda,\mu]].
\]
These two identities make $\Omega(M,\g)$ a Lie superalgebra over $C^\infty(M)$, with grades $\Omega^{\rm even}(M,\g)$ and $\Omega^{\rm odd}(M,\g)$.

An Ehresmann connection $A$ on a principal $G$-bundle is locally a $\g$-valued 1-form, and the exterior covariant derivative is given by:
\[
         d_A \omega = d\omega + [A,\omega].
\]
Obviously the ordinary dif\/ferential $d$, which acts only on form parts:
\[
  d(\omega)= d(\omega^\xa \tensor v_\xa) := d\omega^\xa \tensor v_\xa\, ,
\]
is a graded derivation over the bracket of $\g$-valued forms:
\begin{gather*}
  d[\omega, \mu]    = [d\om, \mu] + (-1)^p[\om, d\mu],
            \qquad
            \om \in \Omega^p(M,\g).
\end{gather*}
Since each $[A,-]$ is also a graded derivation, by the Jacobi identity, so too is $d_A$:
\begin{gather*}
   d_A[\omega, \mu]
            = [d_A\omega,\mu] + (-1)^p[\omega,d_A\mu],
            \qquad
            \om \in \Omega^p(M,\g).
\end{gather*}

Now suppose $\bilin\maps \g\tensor \g \to \R$ is a (nondegenerate) $\Ad(G)$-invariant bilinear form.  As in the usual notation of Yang--Mills theory, $\bilin$ acts on $\g$-valued dif\/ferential forms by
\[
      \bilin(\om\wedge \mu) := \bilin(v_\xa, v_\xb)\, \om^\xa \wedge \mu^\xb.
\]
When $\bilin$ is the Killing form,
this can be written
\[
     \bilin(\om\wedge \mu) = C^\xg_{\xa\xd}C^\xd_{\xb\xg}\, \om^\xa \wedge \mu^\xb,
\]
where $C^\xd_{\xa\xb}$ are the structure constants of $\g$ relative to the basis $v_\xa$.
The invariance of $\bilin$ on $\g$ gives a graded invariance of $\bilin$ for $\g$-valued forms.  In particular, for $\lambda\in \Omega^r(M,\g)$:
\[
    \bilin([\lambda,\om], \mu) = (-1)^{pr+1}\bilin(\om,[\lambda,\mu]).
\]
The exterior derivative is a graded derivation over $\bilin$ in the sense that, for $\om\in \Omega^p(M,\g)$,
\[
     d(\bilin(\om\wedge \mu)) = \bilin(d\om\wedge \mu) + (-1)^p\bilin(\om\wedge d\mu).
\]
But, using the graded invariance under the Lie bracket, we also have
\[
      d(\bilin(\om\wedge \mu)) = \bilin(d_A \om\wedge \mu) + (-1)^p\bilin(\om\wedge d_A \mu),
\]
which lets us integrate by parts using any exterior {\em covariant} derivative.

\section{Special properties of six-dimensional Lie algebras}
\label{apx:Lie-alg}

The 6-dimensional Lie algebras $\so(4)$, $\so(3,1)$, and $\so(2,2)$ inherit a notion of Hodge duality by the fact that they are isomorphic as vector spaces to $\Lambda^2 \R^4$.  This Hodge dual satisf\/ies:
\begin{itemize}\itemsep=0pt
\item $\star [X,Y]=[X,\star Y]$;
\item $\tr(X{\star Y}) = \tr(Y{\star X})$, where `$\tr$' denotes the Killing form;
\item $\star^2 = +1$ for $\g=\so(4)$ or $\so(2,2)$,  $\star^2 = -1$ for $\g=\so(3,1)$.
\end{itemize}
An oft-used consequence of the f\/irst property in calculations in this paper is that if $\om$ is an connection for one of these groups, the covariant dif\/ferential $d_\om$ commutes with the Hodge star:
\[
   d_\om (\star X) = d (\star X) + [\om,\star X]  = \star(dX + [\om,X]) = \star d_\om X.
\]

It is also easy to check from the above properties that $\so(4)$ and $\so(2,2)$ have self dual and anti self dual subalgebras and that these are Killing-orthogonal.  This leads to the Lie algebra coincidences
\begin{gather*}
   \so(2,2)   \iso \Sl(2,\R)\oplus \Sl(2,\R), \qquad
   \so(4)  \iso \su(2)\oplus \su(2).
\end{gather*}

The Hodge star operator extends to $\Iso(2,1)$ and $\Iso(3)$, which are Wigner contractions of the algebras above.  For these, there is a choice, depending on how the contraction is performed, between $\star^2 = \pm 1$.

\begin{lemma}
Let $\g$ be one of the Lie algebras $\so(4)$, $\so(3,1)$, $\so(2,2)$, $\Iso(3)$, or $\Iso(2,1)$.  Then every symmetric invariant bilinear form $\bilin$ on the Lie algebra $\g$ is of the form
\[
        \bilin(X,Y) = \tr(X(c_0 + c_1{\star})Y).
\]
\end{lemma}
\begin{proof}
We have 2-dimensional spaces of invariant inner products, with coordinates $(c_0, c_1)$, on each of these Lie algebras, and adjoint-invariance is a linear condition, so we need only show that the spaces of inner products are 2-dimensional in each case.
The Lie algebras $\so(4)$, $\so(3,1)$ and $\so(2,2)$ correspond to the disconnected Dynkin diagram $D_2$ (they have $\so(4,\C) \iso \so(3,\C) \oplus \so(3,\C)$ as their common complexif\/ication).  Every invariant symmetric bilinear form on one of these Lie algebras is thus a linear combination of Killing forms associated to the two factors.  In particular, the vector space of such invariant bilinear forms on the Lie algebra that are invariant is 2-dimensional.

Next, $\Iso(3)$ is the semidirect sum of $\h = \so(3)$ and the abelian ideal $\p\iso \R^3$,  An invariant inner product on $\Iso(3)$ is thus the sum of invariant inner products on the $\so(3)$ and $\R^3$ parts, and invariant inner products on these are unique up to scale.  The same arguments hold for $\Iso(2,1)$, and in either case we get a 2-dimensional space of invariant inner products.
\end{proof}

Checking which values of $(c_0,c_1)$ result in a degenerate inner product is a straightforward computation.

Starting with a generic nondegenerate invariant symmetric bilinear form on $\g$, we can recover the Killing form as well as the bilinear form of (\ref{CS-star}) with the Hodge star.  For this, we use the symmetric Lie algebra structure of $\g=\h\oplus \p$ and observe that, since $\h$ and $\p$ are Killing-orthogonal,
\[
      \tr(\widetilde X\widetilde Y) = \tr(XY)
\]
while, since $\star$ switches $\h$ and $\p$ components,
\[
     \tr(\widetilde X{\star \widetilde Y}) = - \tr(X{\star Y})
\]
for all $X,Y\in \g$.

\subsection*{Acknowledgements}

I thank James Dolan for many helpful discussions about geometry.  I am also grateful for helpful discussions with John Baez, Steve Carlip, Stanley Deser, Stef\/fen Gielen, Andrew Waldron, and Joshua Willis.
This work was supported in part by the National Science Foundation under grant DMS-0636297.

\pdfbookmark[1]{References}{ref}
\LastPageEnding

\end{document}